\documentclass[12pt]{iopart}
\usepackage{amsthm,amssymb}
\usepackage{graphicx,psfrag,subfigure}
\usepackage{times}

\newif\iffigures
\figurestrue

\newcommand{\Rset}{\mathbb{R}}
\newcommand{\Sset}{\mathbb{S}}
\newcommand{\Tset}{\mathbb{T}}
\newcommand{\Zset}{\mathbb{Z}}

\theoremstyle{plain}
\newtheorem{thm}{Theorem}
\newtheorem{pro}[thm]{Proposition}
\newtheorem{lem}[thm]{Lemma}
\newtheorem{cor}[thm]{Corollary}

\theoremstyle{remark}
\newtheorem{remark}{Remark}

\newcommand{\Constant}{{\rm constant}}

\newcommand{\Diagonal}{\mathop{\rm diag}\nolimits}
\newcommand{\Graph}{\mathop{\rm graph}\nolimits}

\newcommand{\Residue}{{\rm res}}

\newcommand{\am}{\mathop{\rm am}\nolimits}
\newcommand{\sn}{\mathop{\rm sn}\nolimits}
\newcommand{\cn}{\mathop{\rm cn}\nolimits}

\begin{document}

\title[Nonpersistence of resonant caustics in billiards]
      {Nonpersistence of resonant caustics in perturbed elliptic billiards}

\author{Sonia Pinto-de-Carvalho$\dagger$ and Rafael Ram\'{\i}rez-Ros$\ddagger$}

\date{today}

\address{$\dagger$ Departamento de Matem\'atica, ICEx,
                   Universidade Federal de Minas Gerais,
                   CP 702, 30123-970, Belo Horizonte, MG, Brazil}

\address{$\ddagger$ Departament de Matem\`{a}tica Aplicada I,
                    Universitat Polit\`{e}cnica de Catalunya,
                    Diagonal 647, 08028 Barcelona, Spain}

\begin{abstract}
Caustics are curves with the property that a billiard trajectory,
once tangent to it, stays tangent after every reflection at the boundary
of the billiard table. When the billiard table is an ellipse,
any nonsingular billiard trajectory has a caustic, which can be either
a confocal ellipse or a confocal hyperbola.
Resonant caustics ---the ones whose tangent trajectories
are closed polygons--- are destroyed under
generic perturbations of the billiard table.
We prove that none of the resonant elliptical caustics persists under
a large class of explicit perturbations of the original ellipse.
This result follows from a standard Melnikov argument and the analysis of the
complex singularities of certain elliptic functions.
\end{abstract}


\eads{\mailto{sonia@mat.ufmg.br},
      \mailto{Rafael.Ramirez@upc.edu}}

\section{Introduction and main result}

Birkhoff~\cite{Birkhoff1927} introduced the problem of
\emph{convex billiard tables} more than 80 years ago as a way to describe
the motion  of a free particle inside a closed convex smooth curve.
The particle is reflected at the boundary according to the law
``angle of incidence equals angle of reflection''.
Good modern starting points in the literature of the billiard problem
are~\cite{KozlovTreshchev1991,Tabachnikov1995}.

\emph{Caustics} ---curves with the property that a billiard trajectory,
once tangent to it, stays tangent after every reflection---
are the most distinctive geometric objects inside billiard tables,
since they are a geometric manifestation of the regularity of their
tangent trajectories.
For example,
integrable billiards have a continuum of caustics,
whereas the nonexistence of caustics inside a convex billiard table
implies that there are some billiard trajectories
whose past and future behaviours differ dramatically.
See, for instance,~\cite{Mather1982}.
Hence, the existence and persistence of caustics
are two fundamental questions in billiards.
Most of the literature deals with convex caustics,
since they are easier to understand and related to ordered trajectories.
Two exceptions are~\cite[\S 3]{GutkinKatok1995} and~\cite{Knill1998}.

We summarize the classical existence results as follows.
On the one hand, if the boundary curve is smooth enough and strictly convex,
then there exists a collection of smooth convex caustics close to the boundary
of the table whose union has positive area~\cite{Douady1982,Lazutkin1973}.
On the other hand,
Mather~\cite{Mather1982} proved that there are no smooth convex caustics
inside a convex billiard table when its boundary curve 
has some flat point.
Gutkin and Katok~\cite{GutkinKatok1995} gave a quantitative version
of Mather's theorem.

The robustness of a smooth convex caustic is closely related to
the arithmetic properties of its \emph{rotation number},
which measures the number of turns around the caustic per impact.
Caustics with Diophantine rotation numbers persist
under small perturbations of the boundary curve.
This follows from standard KAM arguments~\cite{Douady1982,Lazutkin1973}.
On the contrary, \emph{resonant caustics}
---the ones whose tangent trajectories are closed polygons,
   so that their rotation numbers are rational---
are fragile structures that generically break up.
See, for instance,~\cite{RamirezRos2006}.

This raises two complementary questions.
First,
to characterize the perturbations that preserve/destroy a given
resonant caustic of a billiard table.
Second,
to determine all resonant caustics that are preserved/destroyed
under a given perturbation of an integrable billiard table.
These questions have been studied by several authors.
Baryshnikov and Zharnitsky~\cite{BaryshnikovZharnitsky2006} proved
that the perturbations preserving a given resonant caustic of a
smooth convex billiard table form an infinite-dimensional Hilbert manifold.
As a sample,
we point out that this Hilbert manifold is given by
the set of billiard tables with constant width when the rotation number
of the unperturbed caustic is one half~\cite{Knill1998}.
Concerning the second question,
Ram{\'\i}rez-Ros~\cite{RamirezRos2006} gave a sufficient condition for the
break-up of the resonant circular caustics inside a circular billiard table,
in terms of the Fourier coefficients of the perturbation,
see Remark~\ref{rem:Circles} below.

In this paper we tackle the second question when the billiard boundary
is an ellipse.
In that case,
the billiard dynamics is integrable and
any billiard trajectory has a caustic~\cite{Tabachnikov1995}.
The caustics are the conics confocal to the original ellipse:
confocal ellipses, confocal hyperbolas, and the foci.
Poncelet~\cite{Poncelet1822} showed that if a billiard trajectory
inside an ellipse is a closed polygon,
then all the billiard trajectories sharing its caustic are also
closed polygons.
Even more,
if a billiard trajectory tangent to one of the elliptical caustics
is a \emph{$(m,n)$-gon}
---a closed polygon with $n$ sides that makes $m$ turns around its caustic---,
then all the billiard trajectories sharing its caustic are also $(m,n)$-gons,
and their caustic is called \emph{$(m,n)$-resonant}.
(These two definitions are not restricted to billiards inside ellipses.)
We shall see in Section~\ref{sec:Ellipse} that there is a
unique $(m,n)$-resonant elliptical caustic
for any relatively prime integers $m$ and $n$ such that $1 \le m < n/2$.
Our main result is that all these resonant elliptical caustics break up
under a large class of explicit perturbations of the original ellipse,
see Theorem~\ref{thm:MainTheorem}.

The following notations are required to state the main result.
Once fixed the ellipse
\[
Q =
\left\{
(x,y) \in \Rset^2 : \frac{x^2}{a^2} + \frac{y^2}{b^2} = 1
\right\},\qquad
a > b > 0,
\]
we consider its associated elliptic coordinates $(\mu,\varphi)$
given by the relations
\[
x = c \cosh \mu \cos \varphi,\qquad
y = c \sinh \mu \sin \varphi,
\]
where $c = \sqrt{a^2-b^2}$ is the semifocal distance of $Q$.
The equation of the ellipse $Q$ in this elliptic coordinates
is $\mu \equiv \mu_0$, where $\cosh \mu_0 = a/c$ and $\sinh \mu_0 = b/c$.
Hence, any smooth perturbation $Q_\epsilon$ of the ellipse $Q$
can be written in elliptic coordinates as
\begin{equation}\label{eq:EllipticPerturbation}
\mu = \mu_\epsilon(\varphi) =
\mu_0 + \epsilon \mu_1(\varphi) + \Or(\epsilon^2),
\end{equation}
for some $2\pi$-periodic smooth function $\mu_\epsilon(\varphi)$.

\begin{thm}\label{thm:MainTheorem}
Let $\mu_1(\varphi)$ be a $2\pi$-periodic entire function.
If $\mu_1(\varphi)$ is not constant
(respectively, $\mu'_1(\varphi)$ is not $\pi$-antiperiodic),
then none of the $(m,n)$-resonant elliptical caustics with odd $n$
(respectively, even $n$) persists under the
perturbation~(\ref{eq:EllipticPerturbation}).
\end{thm}

Our proof is based on the study of the persistence of the resonant
rotational invariant circles (resonant RICs) of some twist maps
by means of a first-order Melnikov method.
Only convex caustics can be related to the RICs of those twist maps.
Thus, there is no direct way to extend the same procedure to the
nonconvex caustic hyperbolas,
but we believe that the same results hold for them.

\begin{remark}\label{rem:HigherOrder}
If $\mu_\epsilon(\varphi)$ is constant,
then the perturbed curves $Q_\epsilon$ are ellipses,
so all caustics (resonant or not) are preserved.
Hence, the hypothesis $\mu_1(\varphi)$ nonconstant is natural,
since we are using a first-order method.
Nevertheless, we can still state some results when this hypothesis fails.
More precisely, let us assume that
\[
\mu_\epsilon(\varphi) =
\mu_0 + \epsilon \mu_1 + \cdots + \epsilon^{i-1}\mu_{i-1} +
\epsilon^i \mu_i(\varphi) + \Or(\epsilon^{i+1}),
\]
for some $\mu_0,\ldots,\mu_{i-1} \in \Rset$ and some nonconstant
$2\pi$-periodic entire function $\mu_i(\varphi)$.
Then:
\begin{itemize}
\item
If $n$ is odd,
all the $(m,n)$-resonant elliptical caustics with odd $n$ break up.
This result is a corollary of Theorem~\ref{thm:MainTheorem}.
It suffices to consider $\delta = \epsilon^i$ as the new perturbative parameter,
$Q^\ast_\epsilon =
 \{ \mu \equiv \mu_0 + \cdots + \epsilon^{i-1}\mu_{i-1} \}$
as the unperturbed ellipse,
and to realize that $Q_\epsilon$ is a $\Or(\delta)$-perturbation
of $Q^\ast_\epsilon$ whose first-order term in $\delta$ verifies
the hypotheses of Theorem~\ref{thm:MainTheorem}.
\item
If $n$ is even,
we believe that all $(m,n)$-resonant elliptical caustics also break up,
even if $\mu'_i(\varphi)$ is $\pi$-antiperiodic,
but we should use a second-order Melnikov method in order to prove it.
Unfortunately, the computations become too cumbersome.
\end{itemize}
\end{remark}

\begin{remark}\label{rem:CartesianPerturbation}
If we write the perturbed ellipse $Q_\epsilon$ in Cartesian coordinates as
\[
x^2/a^2 + y^2/b^2 + \epsilon P_1(x,y) + \Or(\epsilon^2) = 1,
\]
then
$2(a^2 \sin^2 \varphi + b^2 \cos^2 \varphi) \mu_1(\varphi) +
 ab P_1(a \cos \varphi, b \sin \varphi) = 0$.
In particular, the function $\mu_1(\varphi)$ is $\pi$-antiperiodic when
$P_1(x,y)$ is odd.
\end{remark}

\begin{remark}\label{rem:Circles}
The case of perturbed circular tables was studied using similar techniques
in~\cite{RamirezRos2006}, but the final result was quite different.
Let us recall it for comparison.
Any billiard trajectory inside a circle of radius $r_0$
has some concentric circle of radius $\sqrt{r^2_0-\lambda^2}$ as caustic,
where $0 < \lambda < r_0$ plays the role of a caustic parameter.
If $\lambda = r_0 \sin (m\pi/n)$, then the circular caustic is $(m,n)$-resonant.
Let us write the perturbed circle in polar coordinates $(r,\theta)$ as
\begin{equation}\label{eq:PolarPerturbation}
r = r_\epsilon(\theta) =
r_0 \big( 1 + \epsilon r_1(\theta) + \Or(\epsilon^2) \big),
\end{equation}
for some smooth function $r_\epsilon: \Tset \to \Rset$.
Let $\sum_{l \in \Zset} \hat{r}_1^l \rme^{\rmi l \theta}$
be the Fourier expansion of $r_1(\theta)$ and $n \ge 2$.
If there exists some $l \in n \Zset \setminus \{0\}$ such that
$\hat{r}_1^l \neq 0$, then the $(m,n)$-resonant circular caustics do not
persist, see~\cite[Theorem 1]{RamirezRos2006}.
In particular,
it is not known if the $(m,n)$-resonant circular caustics with odd
(respectively, even) $n$ break up when $r_1(\theta)$ is not constant
(respectively, $r'_1(\theta)$ is not $\pi$-antiperiodic).
\end{remark}

We complete this introduction with a note on the organization.
In Section~\ref{sec:TwistMaps} we develop a general Melnikov theory
to study the persistence of resonant RICs of twist maps.
The general setup is adapted to billiard maps in Section~\ref{sec:Billiards}.
Finally, Theorem~\ref{thm:MainTheorem} is proved in Section~\ref{sec:Ellipse}
by analysing the complex singularities of certain elliptic functions,
an idea borrowed from~\cite{DelshamsRamirez1996}.

\section{Break-up of resonant invariant curves in twist maps}
\label{sec:TwistMaps}

This section is a generalization of~\cite[\S2]{RamirezRos2006},
although several hypotheses have been weakened.
Namely, the unperturbed map can be nonintegrable,
the resonant invariant circle does not need to be horizontal,
and the shift on the invariant circles can be nonconstant.
In spite of it, the essential idea does not change.
A similar theory is contained in~\cite{Rothos2005}.
For a general background on twist maps we refer to the
book~\cite[\S 9.3]{KatokH1995} or to the review~\cite{Meiss1992}.

Let $\Tset=\Rset/2\pi\Zset$,
and $\pi_1: \Tset \times \Rset \to \Tset$ be the natural projection.
Sometimes it is convenient to work in the universal cover $\Rset$ of $\Tset$.
We will use the coordinates $(x,y)$ for both $\Tset \times \Rset$ and $\Rset^2$.
The lines of the form $x = \Constant$ and $y = \Constant$ will be called
vertical and horizontal, respectively.
A tilde will always denote the lift of a function or set to
the universal cover.
If $g$ is a real-valued function,
$\partial_i g$ denotes the derivative with respect to the $i$th variable.
We will assume that all the considered objects are smooth.
Here, smooth means $C^\infty$.
In particular, all the dependences on the perturbative parameter $\epsilon$
are assumed to be smooth.

We will consider certain diffeomorphisms defined on an open cylinder of
the form $Z = \Tset \times Y$,
for some open bounded interval $Y=(y_-, y_+) \subset \Rset$.
Then $\tilde{Z} = \Rset \times Y$ is an open strip of the plane.
A diffeomorphism $f: Z \to Z$ is called an \emph{area-preserving twist map}
when it preserves area, orientation, and verifies the
\emph{twist condition}
\[
\partial_2 \tilde{\pi}_1 \tilde{f}(x,y) \neq 0,\qquad
\forall (x,y) \in \tilde{Z}.
\]
If the twist is positive (respectively, negative),
then the first iterate of any vertical line tilts to the right
(respectively, left).
We also assume, although it is not essential,
that $f$ verifies some \emph{rigid boundary conditions}.
To be more precise,
we suppose that the twist map $f$ can be extended continuously to the closed
cylinder $\Tset \times [y_-,y_+]$ as a rigid rotation on the boundaries.
That is, there exist some \emph{boundary frequencies} $\omega_\pm \in \Rset$,
$\omega_- < \omega_+$, such that $\tilde{f}(x,y_\pm) = (x+\omega_\pm,y_\pm)$.

Let $D=\{(x,x')\in \Rset^2 : \omega_- < x'-x < \omega_+ \}$.
Then there exists a function $h: D \to \Rset$ such that
$\tilde{f}(x,y)=(x',y')$ if and only if
\begin{equation}\label{eq:ImplicitEquations}
y  =  -\partial_1 h(x,x'), \qquad
y' =   \partial_2 h(x,x').
\end{equation}
The function $h$ is called the \emph{generating function} of $f$.
Besides, if $(x'',y'') = \tilde{f}(x',y')$, then
\begin{equation}\label{eq:DifferenceEquation}
\partial_2 h(x,x') + \partial_1 h(x',x'') = 0.
\end{equation}

We study the dynamics of $f$, but it is often more convenient to work with
the lift $\tilde{f}$, so we will pass between the two without comment and,
in what follows, the lift $\tilde{f}$ remains fixed.

A closed curve $\Upsilon \subset Z$ is said to be a
\emph{rotational invariant circle (RIC)} of $f$ when it is
homotopically nontrivial and $f(\Upsilon) = \Upsilon$.
Birkhoff proved that all RICs are graphs of Lipschitz functions.
See, for instance, \cite[\S IV.C]{Meiss1992}.
Let $\upsilon: \Tset \to Y$ be the Lipschitz function  such that
$\Upsilon = \Graph \upsilon := \{ (x,\upsilon(x)) : x \in \Tset \}$.
If $\upsilon$ is smooth, we say that $\Upsilon$ is a \emph{smooth RIC}.

Twist maps do not form a closed set under composition.
For instance, the square of a twist map is not necessarily a twist map,
and indeed typically it is not.
Nevertheless,
any power of a twist map is \emph{locally twist} on its smooth RICs.

\begin{lem}\label{lem:LocalTwist}
If $\Upsilon = \Graph \upsilon$ is a smooth RIC of
an area-preserving twist map $f: Z \to Z$, then
\[
\partial_2 \tilde{\pi}_1 \tilde{f}^n(x,\tilde{\upsilon}(x)) \neq 0,
\qquad \forall x\in \Rset,\qquad \forall n \ge 1.
\]
\end{lem}
\proof
Given any point $p=(x,\tilde{\upsilon}(x)) \in \tilde{\Upsilon}$,
let $p_j = (x_j,\tilde{\upsilon}(x_j)) = \tilde{f}^j(p)$,
$t_j = (1,\tilde{\upsilon}'(x_j))$, and $v_j = (0,1)$.
We identify the tangent planes $T_p \tilde{Z}$ with
the Euclidean plane $\Rset^2$.
Thus, the vector $t_j$ is tangent to $\tilde{\Upsilon}$ at the point $p_j$
and $v_j$ is a vertical vector at $p_j$.
The linear map $\rmd \tilde{f}^n(p): T_p \tilde{Z} \to T_{p_n} \tilde{Z}$
is the composition of the linear maps
$\rmd \tilde{f}(p_j) : T_{p_j} \tilde{Z} \to T_{p_{j+1}} \tilde{Z}$
for $j=0,\ldots,n-1$.
Let $a_j,b_j,c_j,d_j,\alpha_n,\beta_n,\gamma_n,\delta_n \in \Rset$
be the coefficients such that
\begin{eqnarray*}
\rmd \tilde{f}(p_j): &
t_j \mapsto a_j t_{j+1} + c_j v_{j+1}, &
\quad v_j \mapsto b_j t_{j+1} + d_j v_{j+1} \\
\rmd \tilde{f}^n(p):
& t_0 \mapsto \alpha_n t_n + \gamma_n v_n, &
\quad v_0 \mapsto \beta_n t_n + \delta_n v_n.
\end{eqnarray*}
We note that $b_j = \partial_2 \tilde{\pi}_1 \tilde{f} (p_j)$
and $\beta_n = \partial_2 \tilde{\pi}_1 \tilde{f}^n(p)$.
Let us suppose that the twist is positive, so $b_j > 0$.
We want to prove that $\beta_n > 0$ for any integer $n\ge 1$.
The case of negative twist is completely analogous.

We deduce that $c_j = 0$ from the invariance of $\tilde{\Upsilon}$.
Hence, $\beta_n = \sum_{j=0}^{n-1} D_0^{j-1} b_j A_{j+1}^{n-1}$,
where $D_i^j = \prod_{k=i}^j d_k$ and $A_i^j = \prod_{k=i}^j a_k$.
Besides, we note that $d_j > 0$ because the two components of
$C \setminus \Upsilon$ are invariant.
Finally, we get that $a_j > 0$ from the preservation of orientation.
\qed

Roughly speaking,
a RIC is said to be \emph{resonant} when all its points are periodic,
but we need to be more precise.
Let $(x,y)\in Z$ be a periodic point of the twist map $f$,
and let $n$ be its least period.
Then the exists an integer $m$ such that its lift verifies
$\tilde{f}^n(x,y)=( x + 2\pi m,y)$.
Obviously, $\omega_- < 2\pi m/n < \omega_+$.
Such a periodic point is said to be of \emph{type} $(m,n)$.
A RIC is said to be $(m,n)$-\emph{resonant} when all
its points are periodic of type $(m,n)$.

Let $f$ be an area-preserving twist map with
a $(m,n)$-resonant smooth RIC $\Upsilon = \Graph \upsilon$.
Considering area-preserving twist perturbations of the form
$f_\epsilon = f + \Or(\epsilon)$,
we prove in the following lemma that there exists two graphs
$\Upsilon_\epsilon = \Graph \upsilon_\epsilon$
and $\Upsilon^\ast_\epsilon = \Graph \upsilon^\ast_\epsilon$
$\Or(\epsilon)$-close to $\Upsilon$ and such that
$f^n_\epsilon$ projects the first graph onto the second one along
the vertical direction.

\begin{lem}\label{lem:TwoGraphs}
There exist two smooth functions
$\upsilon_\epsilon,\upsilon^\ast_\epsilon:\Tset \to Y$ defined for
$\epsilon \in (-\epsilon_0,\epsilon_0)$, $\epsilon_0 > 0$, such that:
\begin{enumerate}
\item
$\upsilon_\epsilon(x) = \upsilon(x) + \Or(\epsilon)$ and
$\upsilon^\ast_\epsilon(x) = \upsilon(x) + \Or(\epsilon)$, 
uniformly in $x \in \Tset$; and
\item
$f_\epsilon^n\big(x,\upsilon_\epsilon(x)\big)=
 \big(x,\upsilon^\ast_\epsilon(x)\big)$,
for all $x \in\Tset$.
\end{enumerate}
\end{lem}

\proof
We work with the lift of the maps.
Once fixed an angle $x \in \Rset$,
let $y_0 = \tilde{\upsilon}(x)$ and
\[
\tilde{G}(y,\epsilon) :=
\tilde{\pi}_1 \tilde{f}_\epsilon^n(x,y)  - x - 2\pi m.
\]
This function $\tilde{G}(y,\epsilon)$ verifies the hypotheses of
the Implicit Function Theorem at the point $(y,\epsilon)=(y_0,0)$,
since $\tilde{G}(y_0,0)=0$ and
$\partial_1 \tilde{G}\big(y_0,0\big) =
 \partial_2 \tilde{\pi}_1 \tilde{f}^n(x,\tilde{\upsilon}(x)) \neq 0$,
see Lemma~\ref{lem:LocalTwist}.
Consequently, there exist $\epsilon_0,\eta > 0$ such that
the equation $\tilde{G}(y,\epsilon) = 0$ has exactly one solution
$y_\epsilon = y_0 + \Or(\epsilon)$ in the interval $(y_0-\eta,y_0+\eta)$
for all $\epsilon \in (-\epsilon_0,\epsilon_0)$.
We recall that $\tilde{G}(y,\epsilon)$ had $x \in \Rset$ as an extra parameter,
but it appeared in a $2\pi$-periodic smooth way.
Hence, $\epsilon_0$ and $\eta$ can be taken independent from $x$,
the estimate $|y_\epsilon - y_0| = \Or(\epsilon)$ is uniform in $x$,
and $y_\epsilon$ depends in a $2\pi$-periodic smooth way on $x$.
Finally, set $\tilde{\upsilon}_\epsilon(x) = y_\epsilon$ and then
$\tilde{\upsilon}^\ast_\epsilon(x)$ is determined by means of relation
$\tilde{f}_\epsilon^n\big(x,\tilde{\upsilon}_\epsilon(x)\big)=
 \big(x+2\pi m,\tilde{\upsilon}^\ast_\epsilon(x)\big)$.
The functions
$\tilde{\upsilon}_\epsilon,\tilde{\upsilon}^\ast_\epsilon:\Rset \to Y$
are $2\pi$-periodic and smooth, so they can be projected to two smooth
functions $\upsilon_\epsilon,\upsilon^\ast_\epsilon:\Tset \to Y$
that verify the two claimed properties by construction.
\qed

We say that a $(m,n)$-resonant smooth RIC $\Upsilon$ of a twist map $f$
\emph{persists} under an area-preserving twist perturbation
$f_\epsilon = f + \Or(\epsilon)$ whenever the perturbed map has
a $(m,n)$-resonant RIC $\Upsilon_\epsilon$ for
any small enough $\epsilon$ such that
$\Upsilon_\epsilon = \Upsilon + \Or(\epsilon)$.
The corollary below follows immediately from this definition.

\begin{cor}\label{cor:PersistenceGraphs}
The resonant RIC $\Upsilon$ persists under the perturbation $f_\epsilon$
if and only if $\Upsilon_\epsilon = \Upsilon^\ast_\epsilon$.
\end{cor}

Therefore, it is rather useful to quantify the separation
between the graphs $\Upsilon_\epsilon$ and $\Upsilon^\ast_\epsilon$.

\begin{lem}\label{lem:GraphsPotential}
$\upsilon^\ast_\epsilon(x)- \upsilon_\epsilon(x) =
L'_\epsilon(x)$,
where $L_\epsilon : \Tset \to \Rset$ is a function whose lift is
\begin{equation}\label{eq:CompletePotential}
\tilde{L}_\epsilon(x)=
\sum_{j=0}^{n-1} h_\epsilon
( \bar{x}_j(x;\epsilon),\bar{x}_{j+1}(x;\epsilon)),\qquad
\bar{x}_j(x;\epsilon)=
\tilde{\pi}_1 \tilde{f}^j_\epsilon\big(x,\tilde{\upsilon}_\epsilon(x)),
\end{equation}
and $h_\epsilon$ is the generating function of $f_\epsilon$.
\end{lem}

\proof
As long as confusion is avoided,
we will omit the dependence on $x$ and $\epsilon$.
We introduce the notations
$(\bar{x}_j,\bar{y}_j)=\tilde{f}^j(x,\tilde{\upsilon}(x))$
and $\bar{w}_j = \partial \bar{x}_j/\partial x$ for $j=0,\ldots,n$.
Then $\bar{x}_0 = x$ and $\bar{x}_n = x + 2\pi m$,
so $\bar{w}_0 = \bar{w}_n = 1$.
Besides,
$\bar{y}_0 = \tilde{\upsilon}(x)$ and $\bar{y}_n = \tilde{\upsilon}^\ast(x)$.
From the implicit equations~(\ref{eq:ImplicitEquations}), we get that
$\partial_1 h(\bar{x}_0,\bar{x}_1) = -\bar{y}_0$,
$\partial_2 h(\bar{x}_{n-1},\bar{x}_n)= \bar{y}_n$, and
$\partial_{2} h(\bar{x}_{j-1},\bar{x}_j) +
 \partial_{1} h(\bar{x}_j,\bar{x}_{j+1})=0$
for $j = 1,\ldots,n-1$.
Therefore,
\(
\tilde{L}'(x) =
\partial_1 h(\bar{x}_0,\bar{x}_1) \bar{w}_0 +
\sum_{j=1}^{n-1}
\big(\partial_2 h(\bar{x}_{j-1},\bar{x}_j) +
     \partial_1 h(\bar{x}_j,\bar{x}_{j+1}) \big) \bar{w}_j +
\partial_2 h(\bar{x}_{n-1},\bar{x}_n) \bar{w}_n =
\tilde{\upsilon}^\ast(x) - \tilde{\upsilon}(x)
\).
It is immediate to check that $\tilde{L}:\Rset \to \Rset$ is $2\pi$-periodic,
so it can be projected to a function $L:\Tset \to \Rset$.
\qed

\begin{cor}\label{cor:PersistencePotential}
The resonant RIC $\Upsilon$ persists under the perturbation $f_\epsilon$
if and only if $L'_\epsilon(x) \equiv 0$.
\end{cor}

We shall say that $L_\epsilon : \Tset \to \Rset$ is the
\emph{subharmonic potential} of the resonant RIC $\Upsilon$
under the twist perturbation $f_\epsilon$.
It is rather natural to extract information from the low-order terms
of its expansion
$L_\epsilon(x) = L_0(x) + \epsilon L_1(x) + \Or(\epsilon^2)$.
This is the main idea behind any Melnikov approach to a perturbative problem.
The zero-order term $L_0(x)$ is constant (and so useless), since
$L'_0(x) = \upsilon^\ast_0(x) - \upsilon_0(x) =
           \upsilon(x) - \upsilon(x) \equiv 0$.
We shall say that the first-order term $L_1(x)$ is the
\emph{subharmonic Melnikov potential} of the resonant RIC $\Upsilon$
under the twist perturbation $f_\epsilon$.
The proposition below provides a closed formula for
its computation.

\begin{pro}\label{pro:MelnikovPotential}
If $h_\epsilon = h + \epsilon h_1 + \Or(\epsilon^2)$,
then the lift of $L_1(x)$ is
\[
\tilde{L}_1(x) = \sum_{j=0}^{n-1} h_1(x_j,x_{j+1}),\qquad
x_j = \tilde{\pi}_1 \tilde{f}^j(x,\tilde{\upsilon}(x)).
\]
\end{pro}

\proof
Given any $x\in\Rset$, we set
$x_j = x_j(x) := \bar{x}_j(x;0)$ and
$z_j = z_j(x) := \partial_2 \bar{x}_j(x;0)$ for $j=0,\ldots,n$.
Then the $\Or(\epsilon)$-term of~(\ref{eq:CompletePotential}) is
\begin{eqnarray*}
\fl
\tilde{L}_1(x) & = &
\partial_1 h(x_0,x_1) z_0 +
\sum_{j=1}^{n-1}
\Big( \partial_1 h(x_j,x_{j+1})+ \partial_2 h(x_{j-1},x_j) \Big) z_j +
\partial_2 h(x_{n-1},x_n) z_n + \\
\fl
& & \sum_{j=0}^{n-1} h_1(x_j,x_{j+1}).
\end{eqnarray*}
Using the implicit equations~(\ref{eq:ImplicitEquations}) for the
unperturbed twist map, the first summation vanishes.
The terms $\partial_1 h(x_0,x_1) z_0$ and $\partial_2 h(x_{n-1},x_n) z_n$ also vanish,
since $\bar{x}_0(x;\epsilon) = x$ and
$\bar{x}_n(x;\epsilon) = x + 2\pi m$ for all
$\epsilon \in (-\epsilon_0,\epsilon_0)$.
Besides, $x_j = x_j(x) = \bar{x}_j(x;0) =
\tilde{\pi}_1 \tilde{f}^j(x,\upsilon(x))$.
\qed

The following corollary displays the most important property of the 
subharmonic Melnikov potential in relation with the goals of this paper.

\begin{cor}\label{cor:PersistenceMelnikov}
If $L_1(x)$ is not constant,
then the resonant RIC $\Upsilon$ does not persist under the
perturbation $f_\epsilon$.
\end{cor}

\proof
It follows directly from Corollary~\ref{cor:PersistencePotential}
and the estimate $L_\epsilon = \Constant + \epsilon L_1 + \Or(\epsilon^2)$.
\qed

\section{Break-up of resonant caustics in perturbed billiard tables}
\label{sec:Billiards}

\begin{figure}
\begin{center}
\psfrag{g}{$Q$}
\psfrag{s1}{$\gamma(\varphi)$}
\psfrag{s2}{$\gamma(\varphi')$}
\psfrag{f1}{$\vartheta$}
\psfrag{f2}{$\vartheta$}
\psfrag{f3}{$\vartheta'$}
\psfrag{f4}{$\vartheta'$}
\includegraphics[width=4in]{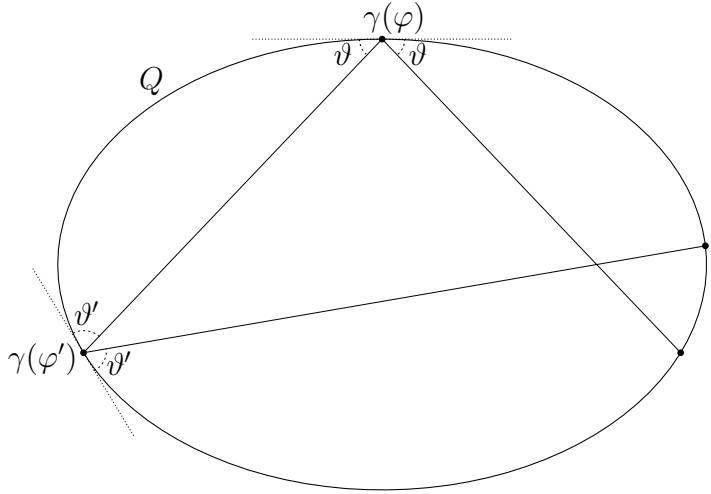}
\end{center}
\caption{The billiard map $f(\varphi,\vartheta) = (\varphi',\vartheta')$.}
\label{fig:BilliardMap}
\end{figure}

Let $Q$ be a closed strictly convex smooth curve in the plane.
Let $\gamma:\Tset \to Q$ be a counterclockwise parametrization.
Let $Z = \Tset \times (0,\pi)$ be an open cylinder.
We can model the billiard dynamics inside $Q$ by means
of a map $f:Z\to Z$, $f(\varphi,\vartheta)=(\varphi',\vartheta')$,
defined as follows.
If the particle hits $Q$ at a point $\gamma(\varphi)$ under
an angle of incidence $\vartheta \in (0,\pi)$ with the tangent vector
at $\gamma(\varphi)$, then, as the motion is free inside $Q$,
the next impact point is $\gamma(\varphi')$, the intersection point with
the boundary and the next angle of incidence is $\vartheta' \in (0,\pi)$,
as in Figure~\ref{fig:BilliardMap}.
A straightforward computation shows that
$f(\varphi,\vartheta)=(\varphi',\vartheta')$ if and only if
\begin{equation}\label{eq:ImplicitEquations2}
|\gamma'(\varphi)| \cos \vartheta  =  -\partial_1 h(\varphi,\varphi'),\qquad
|\gamma'(\varphi')| \cos \vartheta' = \partial_2 h(\varphi,\varphi'),
\end{equation}
where $h: \Tset^2 \setminus \{ \varphi' \neq \varphi \} \to \Rset$ is given by
$h(\varphi,\varphi') = |\gamma(\varphi)-\gamma(\varphi')|$.
Besides, the twist condition holds:
$\partial \varphi'/\partial \vartheta =
 h(\varphi,\varphi')/|\gamma'(\varphi')| \sin \vartheta' > 0$.
Finally, it is geometrically clear that $f$ verifies the
rigid boundary conditions with $\omega_- = 0$ and $\omega_+ = 2\pi$.

A remark is in order. Equations~(\ref{eq:ImplicitEquations2}) differ slightly
from equations~(\ref{eq:ImplicitEquations}),
but identity~(\ref{eq:DifferenceEquation}) still holds and so
the theory developed in the previous section still applies. 

Obviously, one could write the map in the canonical coordinates
---arclength parameter for the boundary and $\cos\vartheta$ as its conjugate---
in order to have $h$ as a generating function,
but this is not a wise choice when dealing with ellipses.

\begin{figure}
\begin{center}
\psfrag{Q}{$Q$}
\psfrag{C}{$C$}
\psfrag{s}{$\gamma(\varphi)$}
\psfrag{fm}{$\vartheta^-$}
\psfrag{fp}{$\vartheta^+$}
\resizebox{3in}{!}{\includegraphics{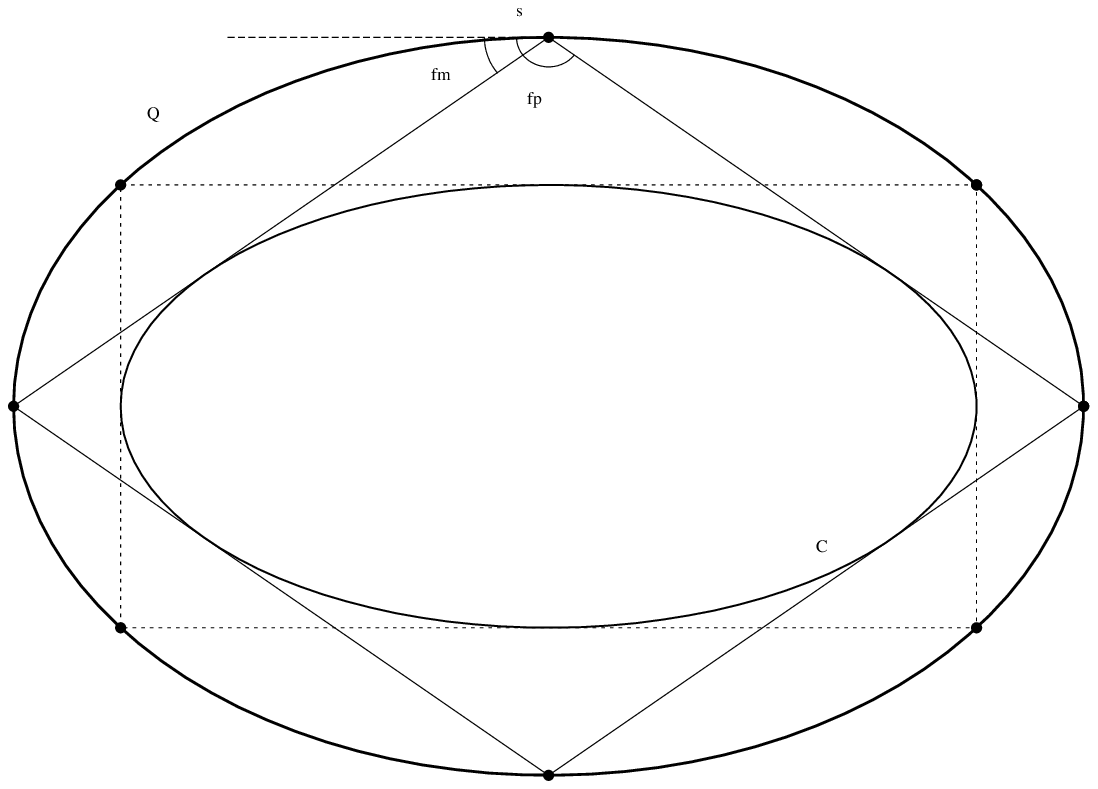}}
\psfrag{0}{$0$}
\psfrag{p/2}{$\frac{\pi}{2}$}
\psfrag{pi/2}{$\pi/2$}
\psfrag{3*pi/2}{$3\pi/2$}
\psfrag{2*pi}{$2\pi$}
\psfrag{pi}{$\pi$}
\psfrag{Yp}{$\Upsilon^+$}
\psfrag{Ym}{$\Upsilon^-$}
\hfill
\resizebox{3in}{!}{\includegraphics{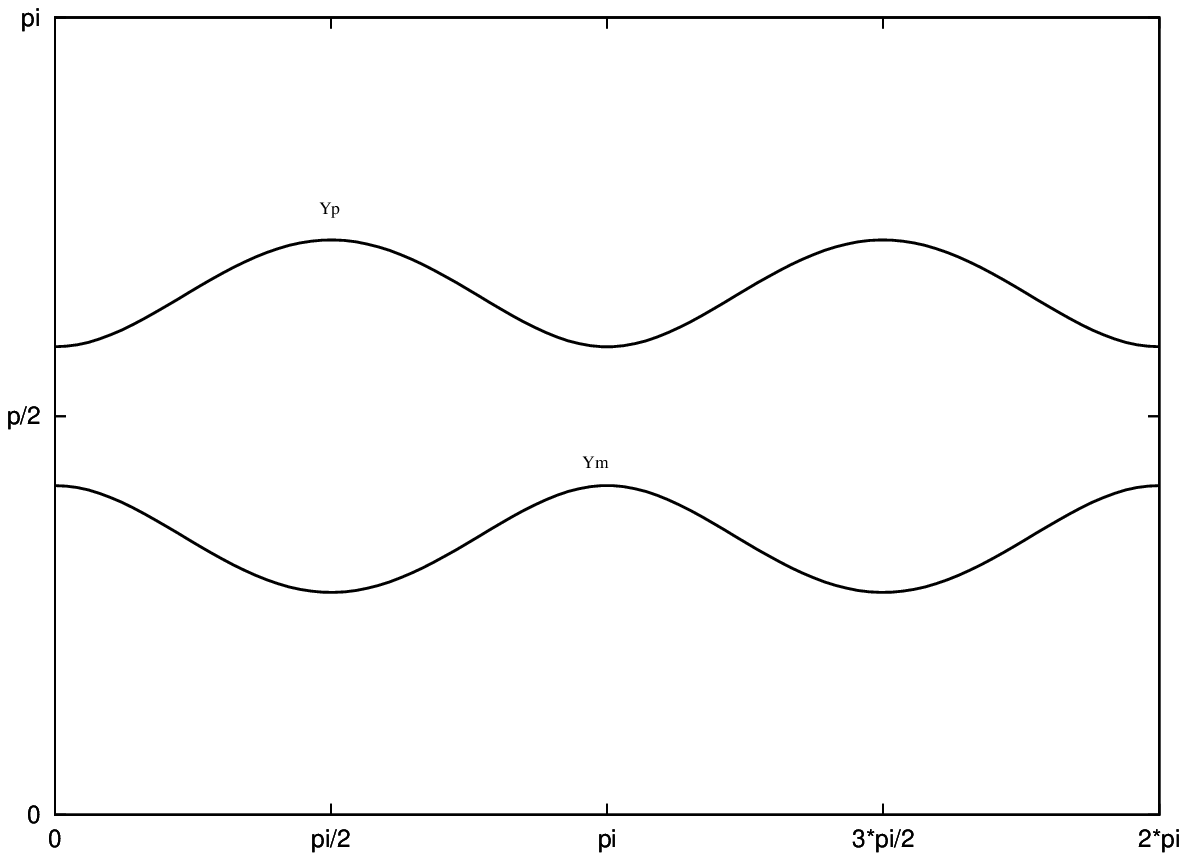}}
\end{center}
\caption{Left: A $(1,4)$-resonant convex smooth caustic $C$.
         Right: Its two smooth RICs $\Upsilon^- = \Graph \vartheta^-$ and
                $\Upsilon^- = \Graph \vartheta^-$ in the phase space
                $Z = \Tset \times (0,\pi)$.}
\label{fig:RICs_Caustics}
\end{figure}

Let us assume that there exists a closed convex smooth caustic $C$
contained in the region enclosed by $Q$.
Then the billiard map $f:Z \to Z$ has two smooth RICs
$\Upsilon^\pm = \Graph \vartheta^\pm \subset Z$.
The functions $\vartheta^\pm: \Tset \to (0,\pi)$ are easy to understand:
$\vartheta^+(\varphi)$ and $\vartheta^-(\varphi)$ are the angles determined
by the two tangent lines to the caustic $C$ from the point
$\gamma(\varphi) \in Q$, see Figure~\ref{fig:RICs_Caustics}.
In particular, $\vartheta^-(\varphi) + \vartheta^+(\varphi) = \pi$.
To fix ideas,
we will assume that $\Upsilon^-$ and $\Upsilon^+$ correspond to the
billiard motion around $C$ in the couterclockwise and clockwise senses,
respectively.
Hence, $0 < \vartheta^-(\varphi) < \pi/2 < \vartheta^+(\varphi) < \pi$.
There is an explicit formula relating the parametrization
of the billiard curve $Q$, the parametrization of the caustic $C$,
and the functions $\vartheta^{\pm}$.
See, for instance, \cite{Douady1982,Knill1998}.

Let $Q$ be a closed strictly convex smooth billiard boundary with a
$(m,n)$-resonant convex caustic $C$,
so that its RIC $\Upsilon^-$ is $(m,n)$-resonant and its RIC $\Upsilon^+$
is $(n-m,n)$-resonant.
We say that $C$ \emph{persists} under a perturbation
$Q_\epsilon = Q + \Or(\epsilon)$ whenever the perturbed billiard curve
has a $(m,n)$-resonant caustic $C_\epsilon$ for any small enough $\epsilon$
such that $C_\epsilon = C + \Or(\epsilon)$.

Let $f_\epsilon$ be the billiard map inside $Q_\epsilon$ and 
$L^-_1(\varphi)$ and $L^+_1(\varphi)$ be the subharmonic Melnikov potentials
of the resonant RICs $\Upsilon^-$ and $\Upsilon^+$
under the area-preserving twist perturbation $f_\epsilon$.
Both potentials coincide, due to the time reversibility of the billiard dynamics.
Therefore, we can skip the $\pm$ signs.
In this context, we will say that $L_1(\varphi)$ is the subharmonic Melnikov
potential of the resonant caustic $C$ for the perturbation $Q_\epsilon$.

\begin{cor}\label{cor:PersistenceMelnikov2}
If $L_1(\varphi)$ is not constant, then the resonant caustic $C$
does not persist under the perturbation $Q_\epsilon$.
\end{cor}

\section{Break-up of resonant caustics in perturbed elliptic billiard tables}
\label{sec:Ellipse}

From now on, we will assume that the unperturbed billiard boundary is
the ellipse
\[
Q =
\left\{
q = (x,y) \in \Rset^2 : \frac{x^2}{a^2} + \frac{y^2}{b^2} = 1
\right\},\qquad
a > b > 0.
\]
It is known that the convex caustics of the billiard inside $Q$ are the
confocal ellipses
\[
C_\lambda =
\left\{
q = (x,y) \in \Rset^2 : \frac{x^2}{a^2-\lambda^2} + \frac{y^2}{b^2-\lambda^2} = 1
\right\},\qquad 0 < \lambda < b.
\]
Let $\rho(\lambda)$ be the rotation number of the elliptical caustic $C_\lambda$.
Then $\rho:(0,b) \to \Rset$ is an analytic increasing
function such that $\rho(0) = 0$ and $\rho(b) = 1/2$.
See, for instance, \cite{CasasRamirez2010}.
Thus, there is a unique $(m,n)$-resonant elliptical caustic for
any relatively prime integers $m$ and $n$ such that $1 \le m < n/2$.
We shall see  that the caustic parameter $\lambda \in (0,b)$ of
the $(m,n)$-resonant caustic is implicitly determined by means of
an equation containing a couple of elliptic integrals,
see equation~(\ref{eq:ResonantCondition}).

The following lemma on elliptic billiards is useful to simplify
the expression of the subharmonic Melnikov potential later on.

\begin{lem}\label{lem:BilliardMotion}
Let $(q_j)_{j \in \Zset}$ be any billiard trajectory inside
the ellipse $Q$ with caustic $C_\lambda$.
Let $p_j = (q_{j+1} - q_j)/|q_{j+1} - q_j|$ be the
unit inward velocities of the trajectory. Then
\[
ab \langle p_{j-1}-p_j, D^{-2}q_j \rangle = 2\lambda,\qquad
\forall j \in \Zset,
\]
where $D=\Diagonal(a,b)$ is the diagonal matrix such that
$Q=\{ q \in \Rset^2 : \langle q,D^{-2} q \rangle = 1\}$.
\end{lem}
\proof
We shall prove that given any point $q = (x,y) \in Q$
and any unit inward vector $p=(u,v) \in \Sset^1$,
the line $\ell = \{ q + \tau p: \tau \in \Rset\}$
is tangent to the conic $C_\lambda$ if and only if
\[
\lambda = - (b xu/a + a yv/b) = -ab \langle p, D^{-2} q \rangle.
\]

To begin with, we note that the line $\ell$ is tangent to the conic $C_\lambda$
if and only if the equation of second order in the variable $\tau$ given by
\[
(x+\tau u)^2/(a^2-\lambda^2) + (y+\tau v)^2/(b^2-\lambda^2) - 1 = 0
\]
has zero discriminant, which is equivalent to the equation
\[
\fl
\left( \frac{xu}{a^2 - \lambda^2} + \frac{yv}{b^2 - \lambda^2} \right)^2 =
\left( \frac{u^2}{a^2 - \lambda^2}  + \frac{v^2}{b^2 - \lambda^2} \right)
\left( \frac{x^2}{a^2 - \lambda^2} + \frac{y^2}{b^2 - \lambda^2} -1 \right).
\]
After some simplifications, we can rewrite this equation as
\[
(xv-yu)^2 =
(b^2 - \lambda^2)u^2 + (a^2 - \lambda^2)v^2 =
a^2 v^2 + b^2 u^2 - \lambda^2,
\]
since $u^2 + v^2 = 1$.
Next, using that $x^2/a^2 + y^2/b^2 = 1$, we obtain that
\[
\lambda^2 = (a^2 v^2 + b^2 u^2)(x^2/a^2 + y^2/b^2) - (xv-yu)^2 =
( b xu/a + a yv/b )^2.
\]
Thus, we have two possibilities:
$\lambda = ab \langle p, D^{-2} q \rangle$ or
$\lambda = -ab \langle p, D^{-2} q \rangle$.
The first one is discarded, because $\lambda > 0$ and
$\langle p, D^{-2} q \rangle < 0$.
The second inequality follows from the fact that
the vector $p$ points inward $Q$ at $q$,
whereas $D^{-2}q$ is an outward normal vector to $Q$ at $q$.

Finally, we note that $-p_{j-1} = (q_{j-1} - q_j)/|q_{j-1} - q_j|$ and
$p_j =(q_{j+1} - q_j)/|q_{j+1} - q_j|$ are the two unit vectors that
point inward $Q$ at the impact point $q_j$ and give the two tangent directions
to the caustic $C_\lambda$.
Therefore,
$\lambda = ab \langle p_{j-1}, D^{-2} q_j \rangle =
-ab \langle p_j, D^{-2} q_j \rangle$.
\qed

\begin{pro}\label{pro:PotentialEllipses}
Let $C_\lambda$ be the $(m,n)$-resonant elliptical caustic confocal
to the ellipse $Q$.
Given any angle $\varphi \in \Tset$,
let $q_j = (a \cos \varphi_j, b \sin \varphi_j)$ be the vertexes of
the $(m,n)$-gon inscribed in $Q$ and circumscribed around $C_\lambda$
such that $q_0 = (a \cos \varphi, b \sin \varphi)$.
Then the subharmonic Melnikov potential of the caustic $C_\lambda$
for the perturbed ellipse~(\ref{eq:EllipticPerturbation}) is
\begin{equation}\label{eq:PotentialEllipses}
L_1(\varphi) = 2 \lambda \sum_{j=0}^{n-1} \mu_1(\varphi_j).
\end{equation}
\end{pro}

\proof
The parametrization of the perturbed ellipse~(\ref{eq:EllipticPerturbation})
is given by
\[
\gamma_\epsilon (\varphi) =
\big( c \cosh \mu_\epsilon(\varphi) \cos \varphi,
c \sinh \mu_\epsilon(\varphi) \sin \varphi \big) =
\gamma_0 (\varphi) + \epsilon \gamma_1(\varphi) + \Or(\epsilon^2),
\]
where $\gamma_0(\varphi) = (a \cos \varphi, b \sin \varphi)$,
$\gamma_1(\varphi) = a b \mu_1(\varphi) D^{-2} \gamma_0(\varphi)$,
and $D = \Diagonal(a,b)$ as above.
The generating function of the billiard map inside the perturbed ellipse is
\[
h_\epsilon(\varphi,\varphi') =
| \gamma_\epsilon(\varphi') - \gamma_\epsilon(\varphi)| =
h_0(\varphi,\varphi') + \epsilon h_1(\varphi,\varphi') + \Or(\epsilon^2).
\]
The first terms of this expansion verify the identities
$h_0(\varphi,\varphi') = |\gamma_0(\varphi') - \gamma_0(\varphi)|$ and
$h_0 (\varphi,\varphi') h_1(\varphi,\varphi') =
\langle \gamma_0(\varphi') - \gamma_0(\varphi),
        \gamma_1(\varphi')-\gamma_1(\varphi) \rangle$.

Let $(q_j)_{j \in \Zset}$ be the billiard trajectory inside the ellipse $Q$
with caustic $C_\lambda$ such that $q_j = \gamma_0(\varphi_j)$
and $\varphi_0 = \varphi$.
The unit inward velocities of this trajectory are
\[
p_j = \frac{q_{j+1} - q_j}{|q_{j+1} - q_j|}
    = \frac{\gamma_0(\varphi_{j+1}) - \gamma_0(\varphi_j)}
           {h_0(\varphi_j,\varphi_{j+1})}.
\]
It follows from Proposition~\ref{pro:MelnikovPotential} that
the subharmonic Melnikov potential is 
\begin{eqnarray*}
L_1 (\varphi) & = & \sum_{j=0}^{n-1} h_1(\varphi_j,\varphi_{j+1}) \\
& = &
\sum_{j=0}^{n-1} \langle p_j,\gamma_1(\varphi_{j+1}) - \gamma_1(\varphi_j) \rangle  \\
& = &
ab \sum_{j=0}^{n-1}
\langle p_j,\mu_1(\varphi_{j+1})D^{-2}q_{j+1} - \mu_1(\varphi_j)D^{-2}q_j \rangle \\
& = &
ab \sum_{j=0}^{n-1}
\langle p_{j-1}-p_j, D^{-2}q_j \rangle \mu_1(\varphi_j) \\
& = &
2 \lambda \sum_{j=0}^{n-1} \mu_1(\varphi_j).
\end{eqnarray*}
We have used the periodicity in the fourth equality and
Lemma~\ref{lem:BilliardMotion} in the last one.
\qed

Next, we give a couple of sufficient conditions
for the subharmonic Melnikov potential to be constant.
These conditions are trivial.
Nevertheless, they play a key role in our problem.
Concretely, we shall check later on that they are also necessary conditions
in the class of $2\pi$-periodic entire functions $\mu_1(\varphi)$.

\begin{cor}\label{cor:Antiperiodic}
Let $\mu_1(\varphi)$ be any $2\pi$-periodic smooth function.
\begin{enumerate}
\item
If the period $n$ is odd, then
$\mu_1(\varphi)$ constant $\Rightarrow L_1(\varphi)$ constant.
\item
If the period $n$ is even, then
$\mu'_1(\varphi)$ $\pi$-antiperiodic $\Rightarrow L_1(\varphi)$ constant.
\end{enumerate}
\end{cor}

\proof
The case $n$ odd is obvious.
If $n$ is even,
the $(m,n)$-gons inscribed in $Q$ and circumscribed around $C_\lambda$
are symmetric with respect to the origin,
so $\varphi_{j+n/2} = \varphi_j + \pi$ and
\[
L'_1(\varphi) = 2\lambda \sum_{j=0}^{n-1} \mu'_1(\varphi_j) =
2\lambda \sum_{j=0}^{n/2 -1}
\left( \mu'_1(\varphi_j) + \mu'_1(\varphi_j + \pi) \right).
\]
In particular, $n$ even and $\mu'_1(\varphi)$ $\pi$-antiperiodic
$\Rightarrow L'_1(\varphi) \equiv 0 \Rightarrow L_1(\varphi)$ constant.
\qed

The subharmonic Melnikov potential of the $(m,n)$-resonant caustic
for the perturbed circle~(\ref{eq:PolarPerturbation}) is
\begin{equation}\label{eq:PotentialCircles}
L_1(\theta) = 2 r_0 \sin(m\pi/n) \sum_{j=0}^{n-1} r_1(\theta_j),\qquad
\theta_j = \theta + 2\pi m j/n,
\end{equation}
see~\cite[Proposition 10]{RamirezRos2006}.
We recall that $\lambda = r_0 \sin (m\pi/n)$ is the $(m,n)$-resonant
caustic parameter of the circle of radius $r_0$.
Besides, all the $(m,n)$-gons inscribed in the
circle of radius $r_0$ and circumscribed around the circle of radius
$\lambda = r_0 \sin (m\pi/n)$ are regular,
so their vertexes are of the form
$q_j = (r_0 \cos \theta_j, r_0 \sin \theta_j)$
with $\theta_j =  \theta + 2\pi m j/n$.
Hence, the function~(\ref{eq:PotentialCircles}) is the limit
of function~(\ref{eq:PotentialEllipses}) when both $a$ and $b$ tend to $r_0$.

Although functions~(\ref{eq:PotentialEllipses}) and~(\ref{eq:PotentialCircles})
look quite similar, they hide a crucial difference.
There is a simple formula for the $\theta_j$ angles,
but not for the $\varphi_j$ ones.
This has to do with the fact that the billiard trajectories inside
a circle of radius $r_0$ sharing a circular caustic with radius
$\lambda = r_0 \sin(\delta/2)$ have a rigid angular dynamics
of the form $\theta \mapsto \theta + \delta$.
On the contrary,
such a rigid angular dynamics does not take place for elliptic tables
when the angle $\varphi$ is considered,
which is a source of technical difficulties in the study of the
subharmonic Melnikov potential~(\ref{eq:PotentialEllipses}).
Nevertheless, it is possible to define a new angular parameter $t$ over
the ellipse $Q$ in such a way that all billiard
trajectories inside $Q$ sharing the elliptical caustic $C_\lambda$
have a rigid angular dynamics of the form $t \mapsto t + \delta$,
for some constant shift $\delta=\delta(\lambda)$.

We need some notations on elliptic functions in order to define
this angular parameter $t$.
We refer to~\cite{AbramowitzS72,WhittakerW27} for a general background
on elliptic functions.
Given a quantity $k\in(0,1)$, called the \emph{modulus},
then $K = K(k) = \int_{0}^{\pi/2}(1-k^2 \sin^2 \phi)^{-1/2}\rmd \phi$
is the \emph{complete elliptic integral of the first kind}.
We also write $K' = K'(k) = K(\sqrt{1-k^2})$.
The \emph{amplitude} function $\varphi = \am t$ is defined through
the inversion of the integral
\[
t = \int_{0}^{\varphi}(1-k^{2} \sin^2 \phi)^{-1/2}\rmd \phi.
\]
Then the \emph{elliptic sinus} and the \emph{elliptic cosinus}
are defined by the trigonometric relations
\[
\sn t = \sin \varphi,\qquad \cn t = \cos \varphi,
\]
respectively.
Dependence on the modulus is denoted by a comma preceding it,
so we can write $\am(t,k)$, $\sn(t,k)$, and $\cn(t,k)$ to avoid any confusion.
In the following lemma it is stated that the angular dynamics becomes
rigid in the angular parameter $t$ given by $\varphi = \am(t,k)$.
It suffices to find the suitable modulus $k$ for
each elliptical caustic $C_\lambda$.

\begin{lem}\label{lem:ChangFriedberg}
Once fixed any caustic parameter $\lambda \in (0,b)$,
we set the modulus $k \in (0,1)$ and the constant shift
$\delta \in (0,2K)$ by the formulae
\begin{equation}\label{eq:ModulusShift}
k^2 = \frac{a^2-b^2}{a^2 - \lambda^2},\qquad
\delta/2 = \int_0^{\vartheta/2} (1-k^{2} \sin^2 \phi)^{-1/2}\rmd \phi,
\end{equation}
where $\vartheta \in (0,\pi)$ is the angle such that
$\sin (\vartheta/2)  = \lambda/b$.
Let
\[
q_j = (a \cos \varphi_j, b \sin \varphi_j) = (a \cn(t_j,k), b\sn(t_j,k))
\]
be any billiard trajectory inside the ellipse $Q$ with caustic $C_\lambda$.
Then $t_{j+1} = t_j + \delta$.
\end{lem}

\proof
By definition, $\varphi_j = \am(t_j,k)$, so
$t_{j+1} - t_{j} =
 \int_{\varphi_j}^{\varphi_{j+1}} (1-k^2 \sin^2 \phi)^{-1/2}\rmd \phi$.
These integrals are equal to a constant $\delta$ that depends only on $C_\lambda$,
see~\cite[page 1543]{ChangFriedberg1988}).
The formula for the constant shift is given in~\cite[page 1540]{ChangFriedberg1988}.
\qed

Remark that if $a = b = r_0$ then the modulus $k$ is equal to zero,
the complete elliptic integral $K$ is equal to $\pi/2$,
the amplitude function is the identity,
the elliptic sinus/cosinus are the usual sinus/cosinus,
the shift $\delta \in (0,\pi)$ is given by $\lambda = r_0 \sin(\delta/2)$,
and the dynamical relation $t_{j+1} = t_j + \delta$
becomes $\varphi_{j+1} = \varphi_j + \delta$.
Thus, we recover the known rigid angular dynamics for circular tables
as a limit of the formulae for elliptic tables.

From now on, $k$ and $\delta$ will denote the modulus and the constant shift
defined in~(\ref{eq:ModulusShift}).
Thus, we shall skip the dependence of the elliptic functions on the modulus.
We note that $C_\lambda$ has eccentricity $k$.
Besides, $C_\lambda$ is the $(m,n)$-resonant elliptical caustic if and only if
\begin{equation}\label{eq:ResonantCondition}
n \delta = 4 K m.
\end{equation}
This identity has the following geometric interpretation.
When a billiard trajectory makes one turn around $C_\lambda$,
the old angular variable $\varphi$ changes by $2\pi$,
so the new angular variable $t$ changes by $4K$.
On the other hand, we have seen that the variable $t$ changes
by $\delta$ when a billiard trajectory bounces once.
Hence, a billiard trajectory inscribed in $Q$ and circumscribed around
$C_\lambda$ makes exactly $m$ turns around $C_\lambda$ after $n$ bounces
if and only if~(\ref{eq:ResonantCondition}) holds.

\begin{pro}\label{pro:Constant}
Let $\mu_1(\varphi)$ be any $2\pi$-periodic entire function.
\begin{enumerate}
\item
If the period $n$ is odd,
then $L_1(\varphi)$ constant $\Leftrightarrow \mu_1(\varphi)$ constant.
\item
If the period $n$ is even, then
$L_1(\varphi)$ constant $\Leftrightarrow \mu'_1(\varphi)$ $\pi$-antiperiodic.
\end{enumerate}
\end{pro}
\proof

Let $\Delta = 2K + 2K'\rmi$ and $z(t) = \cn t + \rmi \sn t$.
If $\varphi = \am t$, then
\begin{eqnarray*}
\rme^{\rmi \varphi} & = &
\cos \varphi + \rmi \sin \varphi = \cn t + \rmi \sn t = z(t), \\
\rme^{-\rmi \varphi} & = &
\cos \varphi - \rmi \sin \varphi = \cn t - \rmi \sn t = z(t + \Delta).
\end{eqnarray*}
We have used that the elliptic cosinus is $\Delta$-periodic,
but the elliptic sinus is $\Delta$-antiperiodic.
We also recall that the elliptic cosinus/sinus are
$2K$-antiperiodic meromorphic functions on the whole complex plane
whose unique singularities are the points of the form
\[
\tau_{r,s} = 2Kr + (1+2s) K' \rmi,\qquad r,s \in \Zset.
\]
Besides, these singularities are just simple poles whose residues are
\[
\Residue(\cn ; \tau_{r,s} ) = (-1)^{r+s+1} \rmi/k,\qquad
\Residue(\sn ; \tau_{r,s} ) = (-1)^{r} /k.
\]
Thus, $z(t)$ is a $2K$-antiperiodic meromorphic function
whose unique singularities are the points of the set
\[
P = \{ \tau_{r,2s+1} : r,s \in \Zset \} =
\tau_\ast + 2K \Zset + 4K'\rmi \Zset,\qquad
\tau_\ast = \tau_{0,-1} = -K'\rmi.
\]
As before, these singularities are just simple poles.

Let $\sum_{l \in \Zset} \hat{\mu}_l \rme^{\rmi l \varphi}$
be the Fourier expansion of $\mu_1(\varphi)$.
Then
\[
\mu_1(\am t) = \mu_1(\varphi) =
\sum_{l \in \Zset} \hat{\mu}_l \rme^{\rmi l \varphi} =
\hat{\mu}_-(z(t+\Delta)) + \hat{\mu}_0 + \hat{\mu}_+(z(t)),
\]
where
$\hat{\mu}_-(z) = \sum_{l=1}^{\infty} \hat{\mu}_{-l} z^l$ and
$\hat{\mu}_+(z) = \sum_{l=1}^{\infty} \hat{\mu}_l z^l$.
We note that the functions $\hat{\mu}_\pm(z)$ are entire,
because $\mu_1(\varphi)$ is entire.
Besides,
\begin{equation}\label{eq:L}
L_1(\am t) = L_1(\varphi) = 2\lambda \sum_{j=0}^{n-1} \mu_1(\varphi_j) =
2\lambda \left( L_-(t) + n \hat{\mu}_0 + L_+(t) \right),
\end{equation}
where
$L_-(t) = \sum_{j=0}^{n-1} \hat{\mu}_-(z(t+\Delta+j\delta))$ and
$L_+(t) = \sum_{j=0}^{n-1} \hat{\mu}_+(z(t+j\delta))$.
Let us study the behaviour of these two functions around the
point $\tau_\ast = -K'\rmi$.
Concretely,
we shall prove that $L_-(t)$ is analytic at $t=\tau_\ast$, whereas
$L_+(t)$ has a nonremovable singularity at $t=\tau_\ast$ provided
$\mu_1(\varphi)$ is nonconstant and $n$ is odd,
or provided $\mu'_1(\varphi)$ is not $\pi$-antiperiodic and $n$ is even.

We begin with a couple of simple observations.
If $j \in \{0,\ldots,n-1\}$, then:
\begin{itemize}
\item[a)]
$\Im (\tau_\ast + \Delta + j\delta) = K'$,
so $\tau_\ast + \Delta + j\delta \not \in P$; and
\item[b)]
$\tau_\ast + j\delta \in P \Leftrightarrow
 4Km j /n = j\delta \in 2K\Zset \Leftrightarrow
 2jm \in n\Zset \Leftrightarrow 2j \in n\Zset \Leftrightarrow
 j \in \{0,n/2\}$.
Here, we have used that $\delta \in \Rset$,
equation~(\ref{eq:ResonantCondition}), and $\gcd(m,n) =1$.
Besides, we stress that the equality $j=n/2$ only can take place
when $n$ is even.
\end{itemize}
We deduce the following results from the above observations.
\begin{itemize}
\item[1)]
$L_-(t)$ is analytic at $t=\tau_\ast$,
because so are $z(t+\Delta + j \delta)$ for $j=0,\ldots,n-1$.
\item[2)]
If $n$ is odd and $\mu_1(\varphi)$ is nonconstant, then:
\begin{itemize}
\item
The function $\hat{\mu}_+(z)$ is nonconstant and entire;
\item 
The function
$L_+(t) - \hat{\mu}_+(z(t)) = \sum_{j=1}^{n-1} \hat{\mu}_+(z(t + j\delta))$
is analytic at $t = \tau_\ast$;
\item
The composition $\hat{\mu}_+(z(t))$ has a nonremovable
singularity at $t=\tau_\ast$; and
\item
The function~(\ref{eq:L}) is nonconstant,
since it has a nonremovable singularity at $t=\tau_\ast$.
\end{itemize}
\item[3)]
If $n$ is even and $\mu'_1(\varphi)$ is not $\pi$-antiperiodic, then:
\begin{itemize}
\item
The sum
$\hat{\sigma}(z) = \hat{\mu}_+(z) + \hat{\mu}_+(-z) =
 2 \sum_{l =1}^{\infty} \hat{\mu}_{2l} z^{2l}$
is a nonconstant entire function;
\item
$z(t + n\delta/2) = z(t+2K m) = (-1)^m z(t) = -z(t)$, since $m$ is odd;
\item
$\hat{\mu}_+(z(t)) + \hat{\mu}_+(z(t+n\delta/2)) = \hat{\sigma}(z(t))$;
\item
The function $L_+(t) - \hat{\sigma}(z(t))$
is analytic at $t = \tau_\ast$;
\item
The composition $\hat{\sigma}(z(t))$
has a nonremovable singularity at $t=\tau_\ast$; and
\item
The function~(\ref{eq:L}) is nonconstant,
since it has a nonremovable singularity at $t=\tau_\ast$.
\end{itemize}
\end{itemize}
Therefore, the proof follows by combining the above results with
Corollary~\ref{cor:Antiperiodic}.
\qed

Finally, we note that our main result
(namely, Theorem~\ref{thm:MainTheorem} stated in the introduction)
follows directly from Corollary~\ref{cor:PersistenceMelnikov2}
and Proposition~\ref{pro:Constant}.

\ack
SP-de-C was partially supported by the CRM and Brazilian agencies
CNPq and FAPEMIG.
RR-R was supported in part by MICINN-FEDER Grant MTM2009-06973 (Spain)
and CUR-DIUE Grant 2009SGR859 (Catalonia).
This work was completed while SP-de-C was a visitor at
the CRM at Barcelona (Spain).
Useful conversations with Pablo S. Casas, Amadeu Delshams, Vadim Kaloshin,
and Vassilios Rothos are gratefully acknowledged.

\end{document}